\documentclass[12pt,a4paper]{amsart}
\usepackage{geometry} 
\usepackage[dvips]{graphicx}
\usepackage{subfigure}
\usepackage{enumerate}
\usepackage{amsmath,amssymb,amsthm}
\usepackage[usenames]{color}

\theoremstyle{plain}
\newtheorem{thm}{Theorem}[section]

\newtheorem{lem}[thm]{Lemma}
\newtheorem{prop}[thm]{Proposition}
\theoremstyle{definition}
\newtheorem{defn}[thm]{Definition}
\theoremstyle{remark}
\newtheorem{rem}[thm]{Remark}

\newtheorem{alg}[thm]{Iterative Step of the Algorithm}

\newtheorem*{ack}{Acknowledgment}

\title{Equal--area method for scalar conservation laws}
\author{Marjeta Kramar Fijav\v{z}, Mitja Lakner, Marjeta \v{S}kapin Rugelj}
\thanks{M. Kramar Fijav\v{z}, M. Lakner, M. \v{S}kapin Rugelj: University of Ljubljana, Faculty of Civil and Geodetic Engineering,
Jamova 2, 1000 Ljubljana, Slovenia.\textit{ E-mails:} mkramar@fgg.uni-lj.si, mlakner@fgg.uni-lj.si, mskapin@fgg.uni-lj.si.}
\date{\today}

\begin{document}
\parindent0cm
\parskip2mm

\begin{abstract}
We study the one-dimensional conservation law. We use a characteristic surface to define a class of functions,
within which the integral version of the conservation law is solved in a simple and direct way.
We develop a simple algorithm for computing the unique solution.
The method uses the equal-area principle and gives the solution for any given time directly.
\medskip

\noindent {\it Key words:} conservation law, equal--area, characteristics, meshfree.\\
\noindent {\it Mathematics Subject Classification (2000):} 35L65, 65M25.

\end{abstract}

\maketitle

\section{Introduction}

We are concerned with an important class of homogeneous hyperbolic differential
equations called \emph{conservation laws.} They state that measurable
quantities do not change in time within an isolated physical system.
We restrict ourselves to the scalar case when the equations are of the form
\begin{equation}\label{de}
u_t + f(u)_x = 0,
\end{equation}
where $u: \mathbb{R}\times[0,\infty)\to\mathbb{R}$ represents the {\em conserved quantity}
while $f:\mathbb{R}\to\mathbb{R}$ is the {\em flux}.  We equip the equation \eqref{de} with the initial condition
 \begin{equation}\label{ic}
 u(x,0)=h(x),\quad x \in\mathbb{R}.
 \end{equation}

In fluid mechanics, equation \eqref{de} with $f(u)=u v$ is called the \emph{equation of continuity} and represents the conservation of mass in the motion of an ideal nonviscous fluid with mass density  $u$ and velocity $v$. For the derivation and many further applications of the equation we refer to \cite{landau}.  Apart from fluid dynamics, this equation is used in various other models for the evolution of continuum quantities such as chemical plug flow reactors or population dynamics.
We shall only mention the Lighthill-Whitham-Richards model \cite{LW,R} which is widely used to describe vehicular traffic flow, see \cite{GP06, Le92, Lo08, whitham}. A very recent application of this model can be found in \cite{data} where the model is compared to experimental data and an algorithm is described to make traffic forecasts.

Let us point out that although the restriction to one dimensional problems may seem oversimplified for practical purposes, many complex problems can be reduced to one-dimensional subproblems. One important example are flows on networks treated in \cite{data,GP06}, some further examples are presented in \cite[Sec. 8]{Se09}. We believe that our approach can be very useful in these problems since it yields a precise solution in any given time.

It is well known that classical (continuously differentiable) solutions of \eqref{de}-\eqref{ic}, even for smooth initial conditions, do not always exist.
In order to allow singularities, which are meaningful for the physical problem behind the equation, one has to generalize  the concept of solutions. The basic idea is to formulate a new extended problem, whose continuously differentiable solutions are the classical solutions to the original equation \eqref{de}.

We first mention the most usual definition of a generalized solution.
A locally $L^1$-function $u$ on $\mathbb{R}\times[0,\infty)$ is called a \emph{weak solution} to  \eqref{de}-\eqref{ic} if
\begin{equation}\label{weak}
\int_{t\geq 0} dt \int_{\mathbb{R}} \left(u \cdot \psi_t + f(u)\cdot \psi_x\right)\; dx +  \int_{\mathbb{R}} h(x)\cdot \psi(x,0)\; dx = 0
\end{equation}
holds for every $C^1$-function $\psi$ on $\mathbb{R}\times[0,\infty)$ with compact support.
Weak solutions are not unique and, in order to obtain the physically correct solution, one has to impose the right \emph{entropy condition}. There is a rich mathematical theory on this topic, see for example monographs \cite{Br05,  Lax, Le92, Le04, Lo08}.

We proceed more directly.
In order to allow discontinuities we consider the {\emph{integral form of the conservation law}
 \begin{equation}\label{if1}
\frac{d}{dt} \int_a^b u(x,t)\; dx = f\left(u(a,t)\right)- f\left(u(b,t)\right)
 \end{equation}
for all $a< b$ and $t>0$ for which $u$ is continuous in points $(a,t)$ and $(b,t)$.
It says that the area  $\int_a^b u(x,t)\, dx$ changes in time
according to the flux at the boundaries. If both $u$
and $f$ are continuously differentiable, (\ref{if1}) implies
(\ref{de}). Mathematically the integral forms \eqref{weak} and \eqref{if1} are equivalent (see ~\cite[p.~28]{Le92}).
We have chosen the latter because of its  direct interpretation   in terms of areas.

We search for the solutions of \eqref{if1} inside a class of functions $\Upsilon$ that is defined by a characteristic surface associated to the problem \eqref{de}-\eqref{ic}. This functions may have jumps along some locally smooth curves in the $(x,t)$-plane (see
Definition \ref{def}).
In Section \ref{Sec_cl} we prove uniqueness of solutions to \eqref{if1} within the class  $\Upsilon$.

Taking the integral over all $\mathbb{R}$, \eqref{if1} implies
\begin{equation}\label{if2}
 \int_{\mathbb{R}} u(x,t)\; dx = \int_{\mathbb{R}} h(x)\; dx
\quad \text{ for all }t\geq 0.
\end{equation}
Hence the area under the graph of the solution does not change in time.
Based on this observation we propose in Section \ref{Sec_method} a simple method  which we call an \emph{equal--area method}. We show  that, under suitable conditions, it yields the unique solution of \eqref{if1} within the  class  $\Upsilon$.

The idea of using the \emph{equal--area principle}  is not new. It has already been exploited in the classical textbook by Whitham \cite[Sect.~2.8-2.9]{whitham}. However, the method there is developed only analytically in a complicated way that is not suitable for explicit computations.  Recently,  Farjoun and  Seibold \cite{Se09} suggested a
\emph{conservative particle method} that uses equal--area
principle. They use the Lagrangean approach, representing the solution
as a cloud of particles which move with the flow. Particles carry
function values and move according to their characteristic
velocities. When the characteristic curves collide, the particles
are merged in such a way that the total area under the function is
conserved.
So far we are not aware of any other numeric method using the equal--area principle.

Let us mention here some classical numerical schemes for conservation laws.
\emph{Finite volume methods}  \cite{Le04} are based on the
integral form (\ref{if1}) instead of the differential equation
(\ref{de}) where the domain is divided into a set of grid cells. The
numerical solution is an approximation of the average value of the
true solution in the grid cell. One of the classical examples of
finite volume methods is the Godunov scheme \cite{God59} that is based
upon the solution  of Riemann problems.

Most numerical methods work well in the interior of the smoothness regions. For
solving problems with discontinuities, \emph{shock tracking} or \emph{front
tracking} \cite{Hol02} methods were developed. They combine some
standard finite difference or finite volume methods in smoothness
regions with an explicit procedure for tracking the location of
discontinuities.

 \emph{Shock--capturing methods} use a different approach
where the goal is to capture discontinuities in the solution
automatically, without explicitly tracking them. One of the modern
shock--capturing methods is a high-resolution flux--limiter Godunov type method
based on solving one--dimensional  Riemann problems which is described in detail in \cite{Le04}.
It is implemented in the software package \emph{Clawpack}  \cite{claw} written by
LeVeque, Berger, et. al., and available on the web.

In Section \ref{Sec_alg} we describe a numerical algorithm we used to implement the equal-area method. We discuss
numerical results obtained by our method  in Section \ref{Sec_numeric} and compare it to the above mentioned finite volume method
\emph{Clawpack}  \cite{claw,Le04} as well as to the newly suggested conservative
particles method \emph{Particleclaw} \cite{particle, Se09}, written by Farjoun and Seibold.

As demonstrated in Section \ref{Sec_numeric}, the presented method performs favorably in terms of efficiency and accuracy.
In contrast to other known methods, where the solution at selected time is obtained by evolution from initial conditions,
it gives the solution for any given time directly. We have to stress, however, that the main purpose of this paper is to
present theoretical principles and results. The  numerical method presented aims to demonstrate the proposed algorithm,
but does not claim necessarily to be the most efficient one.

\section{Characteristic surface and integral solutions}\label{Sec_cl}

 A well-known method for treating the initial value problem (\ref{de})-(\ref{ic}) is
 the {\em method of characteristics}. The characteristics are curves in the $(x,t)$-plane
 along which the function $u$ is constant. In our case they are lines of the form
 \begin{equation}\label{char}
 x_{\xi}(t)=\xi + f'\left(h(\xi)\right)t,\quad \xi\in\mathbb{R},
 \end{equation}
see e.g.~\cite[Sec.~2.2]{Lo08}. It is easy to see that, if a $C^1$-solution $u(x,t)$
of (\ref{de})-(\ref{ic}) exists, its
graph in $\mathbb{R}^3$ is given by
 \begin{equation}\label{graph}
\Gamma:=\{\left(x,t,y \right) = \left( \xi + f'\left(h(\xi)\right)t, t, h(\xi)\right) \mid
\xi\in\mathbb{R}, t\ge 0\}.
 \end{equation}
We shall call $\Gamma$ the {\em characteristic surface} to the problem (\ref{de})-(\ref{ic}). An example
 of such a surface is seen in  Figure \ref{fig:sl:a}.
\begin{figure}[!htbp]
\centering
\subfigure[Parallel characteristics.]
{
\label{fig:sl:a}
\includegraphics[width=65mm]{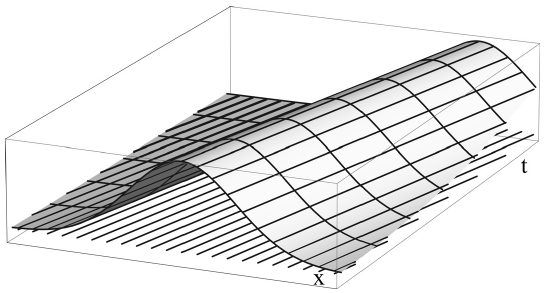}
}
\hspace{3mm}
\subfigure[Characteristics collide.]
{
\label{fig:sl:b}
\includegraphics[width=65mm]{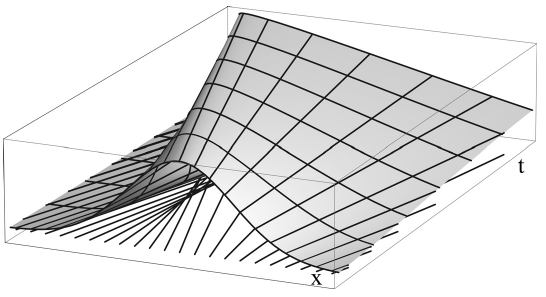}
} \caption{Examples of characteristic surfaces
$\Gamma$.}\label{fig:sl}
\end{figure}

We can form $\Gamma$  a-priori,
before investigating the solvability of \eqref{de}-\eqref{ic}, whereby it
might happen that  it represents a multivalued function which cannot be the
solution of our problem (see Figure \ref{fig:sl:b}). Indeed, this problem occurs
whenever the characteristics collide in the $(x,t)$-plane. In this case the
proper solution has jumps  along
some curves in the $(x,t)$-plane. Its graph, however, is still a
subset of $\Gamma$. Starting with $\Gamma$, the solution can thus
be obtained by finding the appropriate position of the jumps.

In the case when the initial function $h$ is not continuous the above defined characteristic surface $\Gamma$  is not connected. Before proceeding we shall hence modify the definition of characteristic surface to correct this.
For some fixed $t \ge 0$ we define plane transformation $G_t$ as
 \begin{equation}\label{Gt}
G_t(x,y):=(x+f'(y)t,y).
 \end{equation}
Denote by $\gamma_0$ the graph of the initial function $h$ together with vertical line segments joining
discontinuities and
 \begin{equation}\label{curve}
\gamma_t:=G_t(\gamma_0)
 \end{equation}
which is continuous curve for all $t\ge0$, see Figure \ref{fig:claw5}.
\begin{figure}[!htbp]
\centering
\includegraphics[width=140mm]{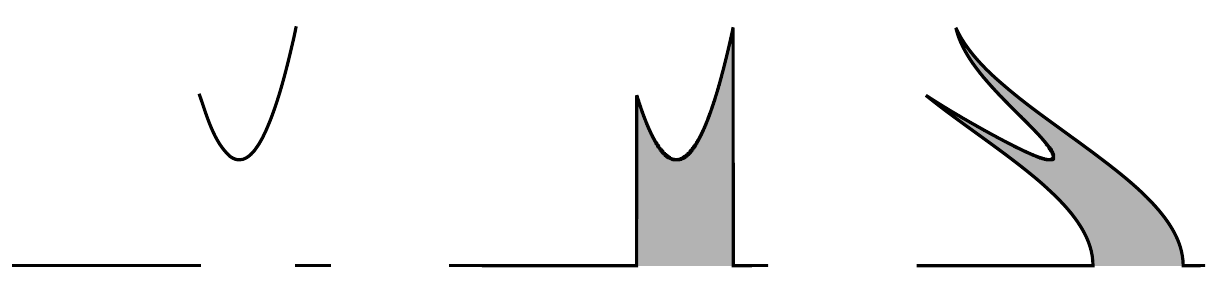}
\caption{The graph of $h$ and the curves $\gamma_0$ and  $\gamma_t$.}\label{fig:claw5}
\end{figure}
\begin{defn}\label{comchar}
The \emph{bounding characteristic surface}  to  the problem \eqref{de}-\eqref{ic} is defined as
\begin{equation}\label{gammahat}
\hat{\Gamma}:=\{\left(x,t,y \right)\mid (x,y)\in\gamma_{t}, t\geq
0\}.
\end{equation}
\end{defn}
Note that $\hat{\Gamma}$ contains the characteristic surface
$\Gamma$ and that the two surfaces agree whenever $h$ is
continuous.

We are now ready to define the appropriate class for our solutions.
\begin{defn}\label{def}
Let $ f\in C^2(\mathbb{R})$ and $h$ is a piecewise $C^1$-function with compact support.
We say that $u\in\Upsilon=\Upsilon(f,h)$  if the following conditions are satisfied:
\begin{enumerate}
\item[(i)] the function $u=u(x,t)$ is defined everywhere on $\mathbb{R}\times[0,\infty)$,
\item[(ii)] the function $u(x,0)=h(x)$ for all $x\in \mathbb{R}$,
\item[(iii)] the graph of $u$ in $\mathbb{R}^3$ is a subset of the bounding characteristic surface $\hat\Gamma$ defined by $f$ and $h$,
\item[(iv)]  the boundary of the graph of $u$ is a finite union of $C^1$-curves, and the projections of these curves on
$(x,t)$-plane intersect any line $t=t_0$ only finitely many times,
\item[(v)] the integral $\int_{-\infty}^{\infty}u(x,t)\; dx$ is a continuous function of $t$ for $t\ge 0$.
\end{enumerate}
\end{defn}

Our aim is to find the solutions to the integral form of the conservation law \eqref{if1} within the  class $\Upsilon$. First we demonstrate that every $u\in\Upsilon$ solves \eqref{de} on its areas of smoothness.

\begin{lem}\label{smooth}
Let $u\in\Upsilon$ be a $C^1$-function on an open set $D\subseteq\mathbb{R}\times[0,\infty)$. Then $u$ is a solution of \eqref{de} on $D$.
\end{lem}
\begin{proof}
If the graph of  $u(D)$ lies on the original characteristic surface $\Gamma$, then $u$ solves \eqref{de} on $D$ by the method of characteristics (see e.g.~\cite[\S 3]{Lax}). Therefore we may assume that the graph of $u(D)$ is contained in $\hat\Gamma\setminus\Gamma$.

We will show that the added vertical lines and their convolutions in time also yield a solution in a similar way as the original characteristics do. For any fixed $\xi$ we parametrize the complemented surface in $\hat\Gamma\setminus\Gamma$ as
$$\left(x,t,u(x,t)\right)=\left(\xi+f'(\tau)t,t,\tau\right). $$
Now by implicit derivation of the $x$-coordinate, by the equality $\tau=\tau(x,t)=u(x,t)$, and by the Implicit function theorem  (see \cite[Theorem 9.28]{rudin}) one obtains
$$u_x=\frac{1}{f''(\tau)t}\quad\text{and}\quad u_t=-\frac{f'(\tau)}{f''(\tau)t},$$
and it is easy to see that $u$, implicitly defined this way, is a local solution of \eqref{de}. \end{proof}

If, however, $u\in\Upsilon$ is not smooth for some $t> 0$ and $x\in\mathbb{R}$, then it has discontinuities called  {\em shocks} which are positioned along piecewise smooth curves $x=s(t)$ in the $(x,t)$-plane, called  {\em shock paths}.
By \ref{def}(iv), there are finitely many shock paths, which are all locally smooth. These paths may have singular points, they may cross or collide, but we can exclude these singularities without loss of generality.

We now continue our treatment by showing that the well-known Rankine-Hugoniot condition holds in our setting.

\begin{lem} \label{solution}
A function $u\in\Upsilon$ is a solution of \eqref{if1}
if and only if the following {\em Rankine-Hugoniot condition}  is satisfied at all the shocks:
\begin{equation}\label{rh}
s'(t) = \frac{f(u^+)-f(u^-)}{u^+-u^-}
\end{equation}
(by $u^+$ and $u^-$ we denote the one--sided limits of the solution $u(x,t)$ from the left and from the right side of the shock, respectively, i.~e.~$u^+(x_0,t)=\lim_{x\searrow x_0} u(x,t)$).
\end{lem}

\begin{proof}
For the sake of simplicity we shall omit the arguments of functions whenever they are clear from the context.
Take any $a < b$. If $u\in\Upsilon$ is smooth for  $x\in(a, b)$ and $t>0$, it solves  \eqref{de} on $(a,b)$ and we have
$$
\frac{d}{dt} \int_a^b u\; dx =  \int_a^b u_t\; dx = - \int_a^b f(u)_x\; dx =f\left(u(a)\right)- f\left(u(b)\right).
$$
Now assume that $u$ has a shock on $(a,b)\times \{t\}$ with shock path  $x=s(t)$. From the condition (iv) of Definition \ref{def} it follows that $x=s(t)$ is a locally $C^1$-function. Hence, we can compute
\begin{align}\label{shock}
\frac{d}{dt} \int_a^b u\; dx &=  \frac{d}{dt} \int_a^{s(t)} u\; dx +\frac{d}{dt} \int_{s(t)}^b u\; dx \notag \\
 &=   \int_a^{s(t)} u_t\; dx +s'(t)u^- +\int_{s(t)} ^b u_t\; dx -s'(t)u^+ \\
 &=  f\left(u(a)\right) - f\left(u^-\right) +  f\left(u^+\right) - f\left(u(b)\right) +
 s'(t)\left(u^--u^+\right).\notag
 \end{align}
 Thus, $u$ is a solution of \eqref{if1} on $(a,b)$ if and only if the Rankine-Hugoniot condition \eqref{rh} is fulfilled at the shock.  If $u$ has more then one shock on $(a,b)\times \{t\}$, we divide the interval according to the shocks and proceed in the same manner. Note that Definition \ref{def}(iv) implies that $u$  has only finitely many shocks.
 \end{proof}

In \cite[Theorem 3.4]{Lax} Lax proved uniqueness of the so-called piecewise generalized solutions to (\ref{de}). We will slightly modify the idea used in his proof in order to show the uniqueness of the solutions in our case.

\begin{prop}\label{unique}
 Let $ f\in C^2(\mathbb{R})$ be strictly concave.  Let $u,v\in\Upsilon$ be two solutions  of  \eqref{if1}
 with the properties
 \begin{enumerate}
 \item[(a)] $u^- < u^+$ and $v^-< v^+$ at every shock for $u$ and $v$, respectively, and
 \item[(b)] $u(x,0)=v(x,0)$ for all $x\in\mathbb{R}$.
 \end{enumerate}
Then the  $L^1$-norm
\begin{equation}\label{l1}
\|u(\cdot,t)-v(\cdot,t)\|_1=0\quad \text{for all }t>0.
\end{equation}
\end{prop}

\begin{proof}
Let $u,v\in\Upsilon$ be  two solutions  to  \eqref{if1} with  $u(x,0)=v(x,0)=h(x)$ for all $x\in\mathbb{R}$. Since they both lie on the same characteristic surface their difference $u-v$ can change  sign only at the shocks of either (or both) of these two functions. Denote these sign-changing curves in $(x,t)$-plane by $y_k(t)$ and order them as
$$y_1(t) <y_2(t)<\cdots <y_{n+1}(t). $$
We work on maximal open intervals of $t$ where the curves $y_k(t)$ do not intersect and the number of these curves does not change.
For almost every $t\geq 0$ we can then write
\begin{equation}\label{pk}
\|u(\cdot,t)-v(\cdot,t)\|_1=\sum_{k=1}^{n} p_k(t) \text{ where }p_k(t)=\int_{y_k(t)}^{y_{k+1}(t)} |u(x,t)-v(x,t)|\; dx.
\end{equation}

Now choose any interval $\left(y_k(t),y_{k+1}(t)\right)$, $1\leq k\leq n$. Without loss of generality we may assume that $u(x,t) > v(x,t)$ on this interval, hence the absolute value under the integral defining $p_k(t)$ can be omitted.

Note that  $p_k(t)$ is  a (continuous) piecewise differentiable function of $t$.  If either $u$ or $v$ has some shocks inside the interval, in each of them the Rankine-Hugoniot condition is satisfied by Lemma \ref{solution}. Dividing the interval according to these shocks and computing the derivative according to this division, similarly as it was done in \eqref{shock}, we see that the values around the shocks cancel out and only the values in  $y_k(t)$ and $y_{k+1}(t)$ are important. Therefore we shall from now on assume that none of $u$ and $v$ has shocks inside the interval $\left(y_k(t),y_{k+1}(t)\right)$.

Thus we may assume that function  $p_k(t)$ is  differentiable and we will now compute its derivative.
As usual we will omit the arguments of functions whenever possible. Since $u$ and $v$ solve  \eqref{de} inside the interval and since the shock paths are piecewise differentiable we have
\begin{align}
p'_k(t) & = 
 \int_{y_k}^{y_{k+1}} (u_t-v_t)\; dx + \left(u^-(y_{k+1})-v^-(y_{k+1})\right)y'_{k+1}  - \left(u^+(y_k)-v^+(y_k)\right)y'_k\notag \\
 & = f\left(u^+(y_k)\right)-f\left(v^+(y_k)\right) -  \left(u^+(y_k) -v^+(y_k)\right)y'_k \label{dpL}\\
 &  -   \left[f\left(u^-(y_{k+1})\right)-f\left(v^-(y_{k+1})\right) - \left(u^-(y_{k+1}) -v^-(y_{k+1})\right)y'_{k+1}\right].\label{dpR}
 \end{align}
From now on we shall explain the calculations only for  the left endpoint of the interval, since the right endpoint can be treated symmetrically. By assumption, $u>v$ inside the interval, hence $u-v$ changes sign at the endpoints and all the jumps are upwards. Taking all this into account we see that  $u$ has a shock in $y_k$ while   $v$ may have a shock or not (in the latter case we take $v^-=v^+$). Furthermore we have
\begin{equation}\label{uv}
u^- \le v^-\le v^+ \le u^+.
\end{equation}
Applying the Rankine-Hugoniot condition for the speed of shock $y'_k$ for $u$ we see that  \eqref{dpL} equals
$$  \left(\frac{f\left(u^+\right)-f\left(v^+\right)}{u^+ -v^+} - \frac{f\left(u^+\right)-f\left(u^-\right)}{ u^+ -u^-} \right) \left(u^+ -v^+\right)\le 0,$$
since by concavity of $f$ and condition \eqref{uv} the first factor is smaller or equal to $0$
while $u^+-v^+>0$.

 Following the same line of arguments for the right endpoint $y_{k+1}$ we obtain the same conclusion for \eqref{dpR}.

We have thus shown  that  $p'_k(t)\leq 0$, $1\leq k\leq n$. By \eqref{pk},  $\|u(\cdot,t)-v(\cdot,t)\|_1$ is then a decreasing function of $t$. Since by assumption  $\|u(\cdot,0)-v(\cdot,0)\|_1=0$  we finally obtain \eqref{l1}.
\end{proof}

We have seen that the solutions  to  \eqref{if1} in the class $\Upsilon$ are unique, provided that the flux is a concave function and the solutions only have jumps upwards. Note that the same is true  for convex flux and downwards jumps.

\begin{rem} Using Lemmas \ref{smooth} and \ref{solution} it is easy to see that our solutions in class $\Upsilon$ are weak solutions (see also \cite[Theorem 4.2]{Br05}).
The condition (a) of Proposition \ref{unique} is actually an entropy condition which yields the physically reasonable solution (see \cite[(3.13)]{Lax}, \cite[p.~36]{Le92}, or \cite[(4.38)]{Br05}).
\end{rem}

\section{The equal--area method}\label{Sec_method}

We now describe the {\em equal--area method} for obtaining the solutions starting from the bounding characteristic surface $\hat\Gamma$ defined in \eqref{gammahat}.  First note that the transformation $G_t$ defined in \eqref{Gt} preserves area, since its Jacobian equals 1.
Hence the area bounded by the curves $\gamma_{t} $ and the $x$-axis remains
unchanged in time and equals the initial area given by $\int_{-\infty}^{\infty}
h(x)\; dx $ (compare the shaded areas in Figure \ref{fig:claw5}).
Intersecting  $\hat{\Gamma}$ with the plane $t=t_0$
for some fixed time $t_0$ yields
 $\gamma_{t_0}$ defined in \eqref{curve}, see also Figure \ref{fig:claw5}.  Our strategy is to insert vertical cuts to  $\gamma_{t_0}$ in
such a way that the areas of the cut-off lobes coincide. Thus the initial area will be preserved.
Carrying out this procedure for all $t$ we obtain a bounded, piecewise continuous function $u(x,t)$, whose graph,
without the added vertical surfaces,  is  contained in ${\hat\Gamma}$, see Figures \ref{fig:claw1} and \ref{fig:claw3}.
\begin{figure}[!htbp]
\centering
\subfigure[The bounding characteristic surface $\hat{\Gamma}$.]
{
\label{fig:claw3:a}
\includegraphics[width=65mm]{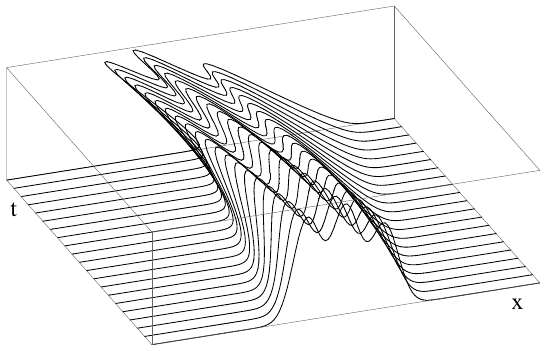}
}
\hspace{5mm}
\subfigure[The solution $u(x,t)$ (with vertical surfaces at the positions of the jumps).]
{
\label{fig:claw3:b}
\includegraphics[width=65mm]{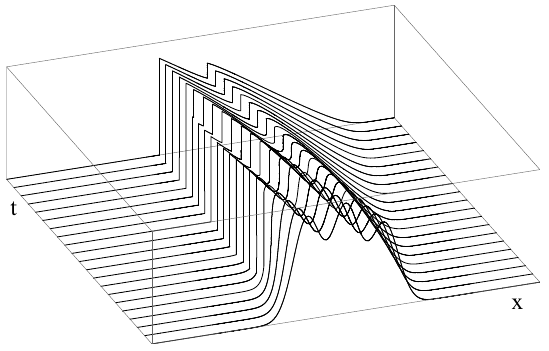}
}
\caption{Graphs for the problem (\ref{de})-(\ref{ic}) for the initial function given in  (\ref{ex_h}). }\label{fig:claw3}
\end{figure}

\begin{lem}\label{implicit}  Let $ f\in C^2(\mathbb{R})$ be strictly concave and let $h$ be a piecewise $C^1$-function with compact support. Then for the solution $u$ obtained by the above described equal--area method the following holds.
\begin{enumerate}
  \item[(a)] At every shock $u$ satisfies the Rankine-Hugoniot condition \eqref{rh} and  has only jumps upwards:
  $u^- < u^+$.
   \item[(b)] The shock paths are piecewise $C^1$-curves.
\end{enumerate}
\end{lem}

\begin{proof}
For $t>0$ we have the curves
$$
\gamma_t(\xi)=(x_t(\xi),y_t(\xi))=(\xi+f'(h(\xi))t,h(\xi)), \qquad \xi\in\mathbb{R}.
$$
For $\xi_1<\xi_2$ we close the ``$S$-curve" $\gamma_t([\xi_1,\xi_2])$ with line segment
between endpoints $\gamma_t(\xi_1)$ and $\gamma_t(\xi_2)$, see Figure \ref{fig:claw11}.
\begin{figure}[!htbp]
\centering
\includegraphics[width=70mm]{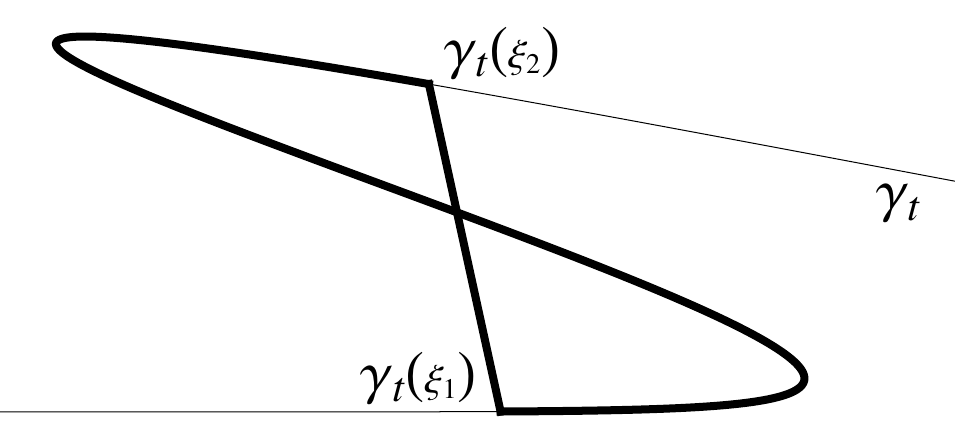}
\caption{
The ``$S$-curve".}\label{fig:claw11}
\end{figure}
The signed area defined by this closed curve by Green's Theorem equals
\begin{align}\label{signed_area}
p_t(\xi_1,\xi_2) &= \frac{1}{2} \int_{\xi_1}^{\xi_2} \left[ x_t(\xi)\;y_t'(\xi)-y_t(\xi)\;x_t'(\xi)\right]d\xi \notag \\
 &+ \frac{1}{2} \int_0^1 \left[ \left(\xi \;x_t(\xi_1)+(1-\xi)\;x_t(\xi_2)\right)\left(y_t(\xi_1)-y_t(\xi_2)\right)\right.\\
 &\left. -\left(\xi \;y_t(\xi_1)+(1-\xi)\;y_t(\xi_2)\right)\left(x_t(\xi_1)-x_t(\xi_2)\right) \right]d\xi.\notag
 \end{align}
We define the mapping $F:\mathbb{R}^2\times\mathbb{R}^+\rightarrow\mathbb{R}^2$ by
$$
F(\xi_1,\xi_2,t):=\left(p_t(\xi_1,\xi_2),x_t(\xi_1)-x_t(\xi_2)\right).
$$
In points where
\begin{equation}\label{enakost}
x_t(\xi_1^0)=x_t(\xi_2^0),
\end{equation}
the determinant of the $2\times2$ Jacobian matrix $\frac{\partial F}{\partial\xi}$ is equal to
$$
2\left(y_t(\xi_2^0)-y_t(\xi_1^0)\right)x_t'(\xi_1^0)x_t'(\xi_2^0)
$$
and is nonzero in jumps given by our method.

If in addition to (\ref{enakost}), we have the equal--area condition $p_{t_0}(\xi_1^0,\xi_2^0)=0$,
then the Implicit function theorem (see \cite[Theorem 9.28]{rudin})
gives the existence of two $C^1$-functions $\widehat{\xi}_1(t)$, $\widehat{\xi}_2(t)$,
such that for all $t$ in some neighborhood of $t_0$ the following holds:
$$
F(\widehat{\xi}_1(t),\widehat{\xi}_2(t),t)=(0,0),\qquad\widehat{\xi}_1(t_0)=\xi_1^0,\qquad \widehat{\xi}_2(t_0)=\xi_2^0.
$$
This means that the equal--area condition holds for all parameters $t$ in this neighborhood.
Moreover, the Implicit function theorem gives us the formula
$$
\left[\begin{array}{cc}\widehat{\xi}'_1(t_0) \\\widehat{\xi}'_2(t_0)\end{array}\right]=
-\left(\frac{\partial F}{\partial\xi}\right)^{-1}\frac{\partial F}{\partial t}.
$$
Using symbolic computation ({\it Mathematica} \cite{mathematica}) one obtains
$$
\widehat{\xi}'_1(t_0)=\frac{f(h(\xi_1^0))-f(h(\xi_2^0))+(h(\xi_2^0)-h(\xi_1^0))f'(h(\xi_1^0))}
{(h(\xi_1^0)-h(\xi_2^0))(1+f''(h(\xi_1^0))h'(\xi_1^0)t_0 )}.
$$
By above, the shock path
$$
x(t)=\widehat{\xi}_1(t)+f'(h(\widehat{\xi}_1(t)))t
$$
is locally a $C^1$-function, which proves (b). Differentiating this function
at the point $t=t_0$, we finally get the Rankine-Hugoniot condition \eqref{rh}
\begin{align}
x'(t_0) &=\widehat{\xi}'_1(t_0)+ f''(h(\xi_1^0))h'(\xi_1^0)\widehat{\xi}'_1(t_0)t_0+f'(h(\xi_1^0)) \notag \\
 &=\frac{f(h(\xi_2^0))-f(h(\xi_1^0))}{h(\xi_2^0)-h(\xi_1^0)}\\
 &=\frac{f(u^+)-f(u^-)}{u^+-u^-}.\notag
 \end{align}
In the case when $f$ is concave,  $f'$ is a decreasing function and the curves $\gamma_t$  are for $t>0$ inclined to the left (regarding $x$-axis). It then  follows from the construction that the solution $u$ only has jumps upwards, hence also the assertion (a) is proved.
\end{proof}

We have proved the local smoothness of the shock paths, a nice property which is often assumed in advance and rarely verified. Assuming a finite number of shocks we are finally  able to prove that the equal-area method yields the unique solution of our problem.

{\begin{thm}  Let $ f\in C^2(\mathbb{R})$ be a strictly concave function and  $h$ a piecewise $C^1$-function  with compact support such that the function
\begin{equation}\label{x}
x_t(\xi)=\xi + f'(h(\xi))t,\quad \xi\in\mathbb{R}
\end{equation}
has only finitely many local extrema for every $t\ge 0$. Then the above described equal--area method yields the unique solution $u\in\Upsilon$ to \eqref{if1}.
\end{thm}}

\begin{proof}
First observe that conditions (i), (ii), (iii), and (v) of Definition \ref{def} are trivially fulfilled for a function $u$ obtained by the equal--area method. Since the number of local extrema of the function \eqref{x} is finite, so is the number of shocks obtained by the equal--area method.  Hence  the condition (iv) of Definition \ref{def} is satisfied by  Lemma \ref{implicit}. Moreover,  using Lemmas \ref{implicit}, \ref{smooth}, \ref{solution}, and \ref{unique} we see that $u\in\Upsilon$ is the unique solution of \eqref{if1}.
\end{proof}

A brief comment is in order here. For our theoretical approach we need to assume that there are only finitely many shocks. We are not aware of any explicit conditions  in terms of functions $f$ and $h$ to meet this assumption (some generic conditions are given in \cite{golubitsky, schaeffer}).
In practice however, it is not difficult to check the finiteness of the number of local extrema of the function  $x_t(\xi)$ in \eqref{x} for any given $f$, $h$, and $t$ (note that our procedure gives a solution for  any fixed time $t$!). Moreover, numerically this condition is always satisfied, since we use  polygonal approximation for  continuous curves.

\section{The algorithm}\label{Sec_alg}

We shall now describe the algorithm based on the procedure
introduced in the previous section. The solution at some
fixed time $t_0$ is obtained directly, without the need to march forward in
time.

We start by taking a polygonal approximation $K_0$ for the
continuous curve
 \begin{equation}
\gamma_{t_0}:=G_{t_0}(\gamma_0),
 \end{equation}
see (\ref{Gt})-(\ref{curve}). We used points on the curve $\gamma_{t_0}$ obtained equidistant parameters on x-axis:
$G_{t_0}(x_i, h(x_i))$. Traveling along the curve we gradually `equalize' the areas.
The obtained graph of solution is a subset of $K_0$ (see Figure \ref{fig:claw1}).

\begin{alg}\label{alg} Let the parametrization $(x(\tau),y(\tau))$, $\tau\in\mathbb{R}$,
of the curve $K_i$ on the $i$-th step be such that $x(\tau)$ is increasing on the far ends of
the interval $(-\infty,\infty)$. First we define three significant  points for the curve $K_i$ (see Figure \ref{fig:claw0}):
\begin{enumerate}
\item  $\beta= x(\tau_1)$ is the first local maximum of  $x(\tau)$ in the direction of the increasing parameter $\tau$.
If such a maximum does not exist, we are done and $K_i$ is the graph of the weak solution  (with redundant vertical lines in the jumps).
\item   $\alpha= x(\tau_2)$ is  is the first local minimum of  $x(\tau)$ from $\tau_1$ onwards.
Since the function  $x(\tau)$ in not bounded from above, such a minimum always exists.
\item $\gamma$ is the minimum of $\beta$ and the first next local maximum of  $x(\tau)$. If such a maximum does not exist, let $\gamma=\beta$.
\end{enumerate}
\begin{figure}[!htbp]
\centering
\includegraphics[width=100mm]{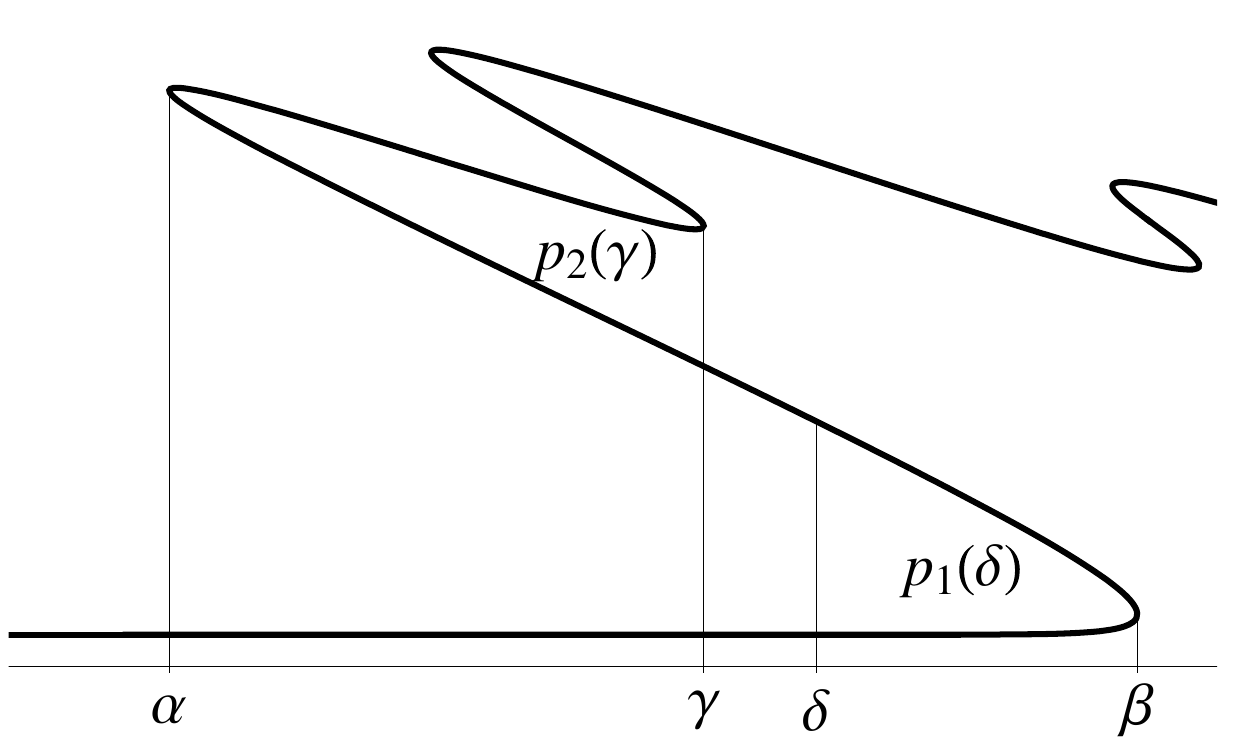}
\caption{Significant points $\alpha$, $\beta$, $\gamma$, and $\delta$ on the curve $K_i$ and the areas $p_1(\delta)$
and $p_2(\gamma)$.}\label{fig:claw0}
\end{figure}
Now we compute the areas. For any $x_0\in(\alpha,\beta)$  denote by $p_1(x_0)$ the area bounded by the line $x=x_0$
and the part of the curve $K_i$ that contains $(x(\tau_1),y(\tau_1))$.
For any $x_0\in(\alpha,\gamma)$  denote by $p_2(x_0)$ the area bounded by the line $x=x_0$ and the part of the
curve $K_i$ that contains $(x(\tau_2),y(\tau_2))$. For  $x\in(\alpha,\gamma)$  let
\begin{equation}\label{p}
 p(x):= p_2(x)-p_1(x).
\end{equation}
Then $p_1(x)$ is continuously decreasing while $p_2(x)$ and $p(x)$
are continuously increasing functions and $p(\alpha)=0-p_1(\alpha)
< 0$. We distinguish two cases:
\begin{enumerate}
\item If $p(\gamma)\geq 0$, let $\delta\in(\alpha,\gamma]$ be the only zero of the function $p(x)$,
therefore $p_1(\delta)=p_2(\delta)$. The curve $K_{i+1}$ is obtained from $K_i$ where the parts of $K_i$
that determine $p_1(\delta)$ and $p_2(\delta)$ are replaced by the vertical line.
\item  If $p(\gamma)< 0$, then $p_1(\gamma)> p_2(\gamma)$ and by continuity and monotonicity  of the
function $p_1(x)$ there exists only one  $\delta\in(\gamma,\beta)$ which satisfies $p_1(\delta)=p_2(\gamma)$
(note that $p_1(\beta)=0$). The point $\delta$ together with the areas  $p_1(\delta)$ and $p_2(\gamma)$ is
marked on the Figure \ref{fig:claw0}. The new curve $K_{i+1}$ is obtained from $K_i$ by replacing those parts
of $K_i$ that determine $p_1(\delta)$ and $ p_2(\gamma)$ by a vertical line.
\end{enumerate}
\end{alg}

In Figure \ref{fig:claw1} there is an example of the resulting
steps of the above algorithm applied to the function given in
(\ref{ex_h}) in time $t_0=4.25$.
\begin{figure}[!htbp]
\centering
\includegraphics[width=120mm]{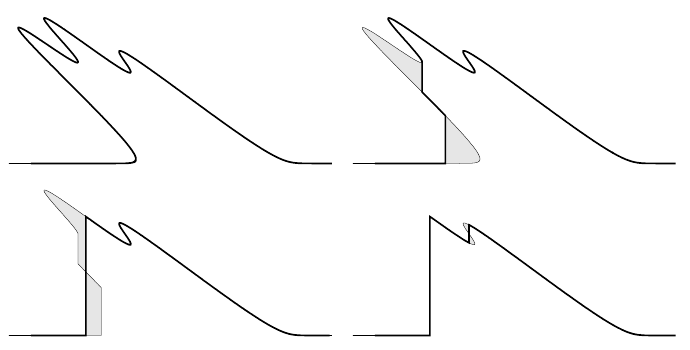}
\caption{Curves $K_i$, $i=0,1,2,3$, obtained in three consecutive steps of Algorithm \ref{alg} resulting
in the  solution at time $t_0$.}\label{fig:claw1}
\end{figure}

We shall briefly describe the method we use to compute the areas needed on each step of Algorithm \ref{alg}.
Let $D$ be a polygon, determined by the points $T_1(x_1,y_1),\dots,T_n(x_n,y_n)$ (we orient them in
counterclockwise direction) and let $T_0(x_0,y_0)$ be any point in the plane. Then the signed area of the triangle $T_0T_iT_{i+1}$ can be computed by
\begin{equation}\label{p012}
p_{0,i,i+1} := \frac{1}{2}\left|\begin{array}{cc}x_i-x_0 & y_i-y_0 \\x_{i+1}-x_0 & y_{i+1}-y_0\end{array}\right|.
\end{equation}
By Green's Theorem one can easily see that the area of the polygon $D$ then equals
\begin{equation}\label{poly}
p=\sum_{i=1}^n p_{0,i,i+1}
\end{equation}
where $T_{n+1}=T_1$.

\section{Numerical results}\label{Sec_numeric}

We have programmed our equal--area method in \emph{Mathematica
}\cite{mathematica} and first compared the results with the basic
Godunov method (which we have also  implemented in
\emph{Mathematica}). Further we have compared our method to an
advanced Godunov method, which is a basis of the widely-used
software package \emph{Clawpack}  \cite{claw}. Finally, we have
made a comparison to the very recent software package
\emph{Particleclaw} \cite{particle}, which uses a Langrangean
particle method and some information on the characteristics. Both
software packages are freely available on the web.

The results of these tests are very good. Our algorithm performs favorably both in terms of time efficiency and accuracy.
Figure \ref{fig:claw2} contains graphs of the solution obtained
by our method and both above mentioned software packages for the initial condition
\begin{equation}\label{ex_h}
 h(x)=\left\{\begin{array}{ll}0.9 e^{-x^2} + 0.7 e^{-(x - 2)^2} + 0.85 e^{-(x + 2)^2} , & x\in [-10,10],\\
 0,& \text{otherwise}.\end{array}\right.
\end{equation}
and the flux function $f(u)=u(1-u)$.

\begin{figure}[!htbp]
\centering
\subfigure[Clawpack]
{
\label{fig:claw2:a}
\includegraphics[width=65mm]{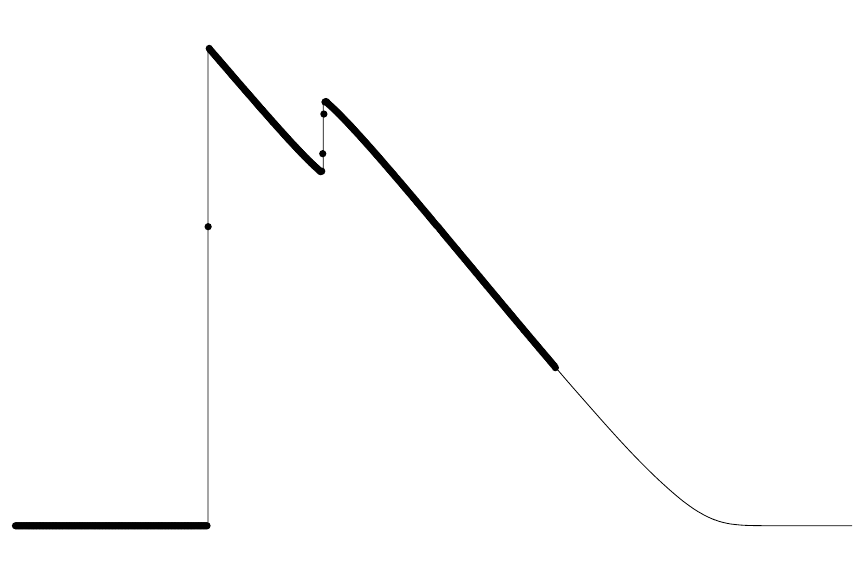}
}
\hspace{7mm}
\subfigure[Particleclaw]
{
\label{fig:claw2:b}
\includegraphics[width=65mm]{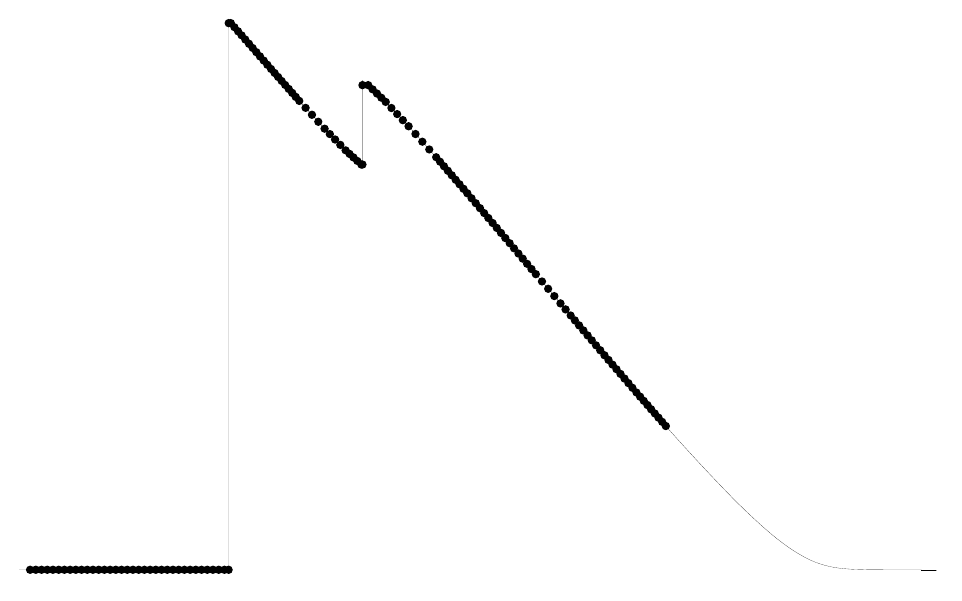}
}
\caption{The solution obtained by the equal--area method (thin line) compared to  \emph{Clawpack}  and \emph{Particleclaw}, respectively  (thick dots). }
\label{fig:claw2}
\end{figure}

Time complexity of our algorithm depends mostly on finding zeros of a function, obtained by the computation of areas of polygons.
We used the secant method and typically $7$ to $12$ iterations (9 on average) were needed
for $10^{-14}$ accuracy. The number of necessary steps of the Algorithm \ref{alg} is bounded above by the
number of stationary points of the function $x_{t_0}(\xi)$ defined in \eqref{x}.

We can approximate the error of the position of the shock. Assuming that the original curve $\gamma_{t_0}$ lies in an
$\varepsilon$-neighborhood of polygonal approximation line,
 we have an approximation of the area
between the ``$S$-curves" as in Figure \ref{fig:claw0}:
$$
l\; \varepsilon\doteq\Delta x\; s
$$
where $l$ is the length of the ``$S$-curve", $\Delta x$ is the displacement of the true shock and $s$ is the height
of the shock. This gives an approximation of the displacement $\Delta x$:
$$
\Delta x \doteq \varepsilon \frac{l}{s}.
$$
In Figure \ref{fig:claw1} polygonal approximation with 1000 points has $\varepsilon$ less than $6.10^{-4}$ and $\Delta x$
is approximated with $3.10^{-3}$.
The method is quadratical, i.e.~doubling the number of points would result in decreasing the value of $\varepsilon$ by factor $4$.

\section{Conclusions}

Using (bounding) characteristic surface we have defined the proper solution to the integral form of the conservation law \eqref{if1}. We have proposed an equal-area method and shown that the obtained (unique) solution has all the desired properties:
it solves \eqref{de}-\eqref{ic} exactly wherever it is smooth, it satisfies the Rankine-Hugoniot condition \eqref{rh} in all the shocks,
and the shock paths are locally smooth. Finally, we have described an algorithm for implementing our method and compared it to some other known methods.

We see the following advantages of our equal--area method.
\begin{itemize}
\item Contrary to classical numerical schemes for conservation laws, it is mesh-free.
\item The method is by its nature exactly conservative.
\item  The solution is computed for any given fixed time. Hence, the errors do not accumulate in time.
\item The method is accurate -- the quality of the approximation relies only on the quality of the starting
approximation of the curve $\gamma_{t_0}$ with a polygonal line
$K_0$.
\item Some methods treat the rarefaction waves separately, using different techniques than for the case of the shock waves. There is no need for  that in our case, and the rarefaction waves are created on the way where appropriate.
\item The obtained shocks are sharp and propagate with correct speed. Their position is obtained automatically by equalizing the
appropriate areas.  Moreover, the shock paths are obtained easily
 by computing the solution for some selected times and then simply projecting  the  shocks from
 the surface to the $(x,t)$-plane (see Figure \ref{fig:claw8}).
\end{itemize}

\begin{figure}[!htbp]
\centering
\includegraphics[width=30mm]{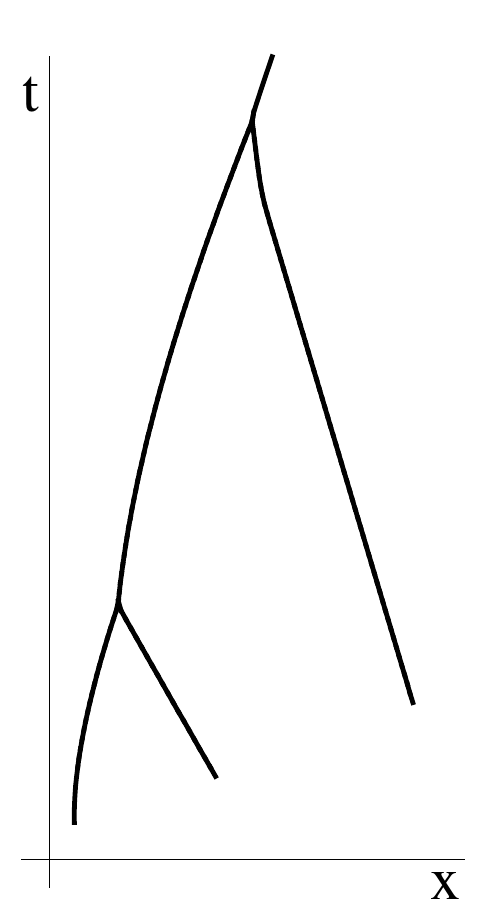}
\caption{Shock paths in $(x,t)$-plane for the initial function given in (\ref{ex_h}).}\label{fig:claw8}
\end{figure}

\begin{ack}
The authors would like to express their gratitude to Ale\v{s} Zalo\v{z}nik and Marijan \v{Z}ura for motivation and many helpful discussions.
\end{ack}

\end{document}